\documentclass[11pt]{article}
\usepackage{amsmath,amssymb,theorem}

\newcommand{\tmop}[1]{\ensuremath{\operatorname{#1}}}

\newcommand{\tmtextit}[1]{{\itshape{#1}}}

\newenvironment{proof}{\noindent\textbf{Proof\ }}{\hspace*{\fill}$\Box$\medskip}

\newtheorem{theorem}{Theorem}[section]
\newtheorem{proposition}[theorem]{Proposition}
\newtheorem{corollary}[theorem]{Corollary}
\newtheorem{definition}[theorem]{Definition}
\newtheorem{lemma}[theorem]{Lemma}
\newtheorem{remark}[theorem]{Remark}
\newtheorem{example}[theorem]{Example}

\def \d {{\mathrm {d}}}

\begin{document}

\begin{center}
{\large{\bf Dirac pairs}}

{\medskip}

Yvette Kosmann-Schwarzbach

{\bigskip}

{\small Centre de Math\'ematiques Laurent Schwartz\\
\'Ecole Polytechnique\\ 
91128 Palaiseau, France}

\texttt{yks@math.polytechnique.fr}

{\bigskip}

{\it Dedicated to Tudor Ratiu}

\end{center}

\bigskip

\begin{abstract} We 
extend the definition of the Nijenhuis torsion of an endomorphism
of a Lie algebroid to that of a relation, 
and we prove that the torsion of the 
relation defined by a bi-Hamiltonian structure vanishes. 
Following Gelfand and Dorfman, we then define Dirac pairs, and we
analyze the relationship of this
general notion with the various kinds of compatible structures 
on manifolds, more generally on Lie algebroids.
\end{abstract}

\section*{Introduction}

What came to be known as `Poisson geometry'
has been essential to the
understanding of mechanics for a very long time: Poisson parentheses,
Poisson brackets, Poisson structures, Hamiltonian structures and Poisson
bivectors are terms familiar from the work of 
Poisson himself, Jacobi, Hamilton and,
in the twentieth century, Tulczyjew, Kirillov and Lichnerowicz.
The introduction and use of 
bi-Hamiltonian structures and Nijenhuis operators in the study of
the integrable systems of mechanics and field theory, due to Magri,
Gelfand and Dorfman, and Fokas and Fuchssteiner, in the late 1970's,
has become the subject of far too many publications to be cited here.

More recently, there appeared what can be called `Dirac geometry', 
in the work of Dorfman
in an algebraic framework suitable for infinite-dimensional dynamical
systems \cite{D1987}, and that of Courant and Weinstein \cite{Courant}
in the geometry of what was later to be 
called the `generalized tangent bundle' of a
smooth manifold. Many advances and generalizations
have appeared, extending the definition and study of Dirac structures to Lie
algebroids, then to Lie bialgebroids \cite{LWX} \cite{Liu} and eventually 
to proto-bialgebroids \cite{YinHe}.
Dorfman, in her 1987 article \cite{D1987},
introduced the Dirac pairs that generalized the bi-Hamiltonian
structures, and she developed their applications in 
her book on the integrability of
nonlinear evolution equations \cite{D}. The purpose of the present
paper is to expand her
approach by introducing Nijenhuis relations in Lie or, more generally, Leibniz
algebras, in order to define Dirac pairs in the double of a Lie
algebroid, and to analyze the relationship of this notion
with the various compatible structures on manifolds and, more
generally on Lie algebroids, that have been the subject of many
papers, most recently~\cite{KSR}.

We have tried to follow Dorfman's terminology as closely as possible, and
yet we were obliged to deviate
from it when necessary for the
clarity of our exposition.
For instance, we have reserved the term `symplectic' for what she called
`invertible symplectic'.

In Section \ref{sec1}, we review definitions 
concerning set-theoretic relations, and we define the Nijenhuis
relations by the vanishing of their torsion.
Section \ref{Hpairs} contains a brief review of bi-Hamiltonian
structures -- which we call `Hamiltonian pairs'--, 
the definition of the relation
defined by two Poisson bivectors and the computation of the torsion
of that relation (Theorem \ref{graphs}), 
which yields the proof that the relation defined by a
Hamiltonian pair is a Nijenhuis relation. 
We define Poisson pairs by the requirement
that a pair of Poisson bivectors define a Nijenhuis relation (Definition
\ref{Ppair}) and, in Corollary \ref{Ppairbis}, 
we compare Hamiltonian pairs and Poisson pairs.
In Section \ref{3}, following Dorfman, 
we introduce the compatibility condition that
defines Dirac pairs (Definition \ref{Dpair})
in such a way that Poisson 
structures constitute a Poisson pair if and only if their graphs
constitute a Dirac pair.
(See Remark \ref{remark} for a comparison with \cite{D} and the rebuttal of
a claim made there concerning Hamiltonian pairs.)
The general definition of Dirac pairs permits the introduction of the
notion of a presymplectic pair. As
expected, in the non-degenerate case, Poisson pairs are in one-to-one
correspondence with presymplectic pairs (Theorem \ref{1-1}). 
In Section \ref{4}, we consider pairs of a Poisson structure together with a
presymplectic structure (P$\Omega$-structures) and pairs of a
presymplectic structure and a Nijenhuis structure ($\Omega$N-structures),
and we charaterize them in terms of Dirac pairs. In Theorem
\ref{oldtheorem} we show that any P$\Omega$-structure defines a Dirac
pair, while in Theorem 
\ref{newtheorem} we prove that for a $2$-form $\omega$ and a Nijenhuis
tensor $N$ in a class of $\Omega$N-structures
containing the symplectic-Nijenhuis structures, the graphs of $\omega$
and $\omega \circ N$ constitute a Dirac pair.

In Section \ref{appendix}, we have included an Appendix in three parts. 
Section \ref{A.1} presents a review of elements of the theory of Lie
algebroids and bialgebroids in order to make the paper essentially 
self-contained. Section \ref{A.2} is a brief review of Terashima's
Poisson and presymplectic functions which we use in Section \ref{A.3} 
to give a new formulation (Theorem \ref{proppoisson})
and a slightly more conceptual proof of a theorem of Yin and He
\cite{YinHe} that characterizes the Dirac structures defined by
bivectors, and we state the dual result that characterizes the Dirac structures
defined by $2$-forms (Theorem \ref{propsympl}).  

Including
Dirac-Nijenhuis manifolds \cite{H} in the framework of Dirac pairs
should be the object of future work, as well as 
extending the definition and study 
of Dirac pairs in the case of the double of a Lie bialgebroid and, 
more generally still, in the case of an artbitrary Courant
algebroid, and establishing the link with the Nijenhuis
structures studied in \cite{C} \cite{ykspreprint}, in the hope of applying these
abstract notions to concrete problems of integrability.

\section{Relations}\label{sec1}

\subsection{Composition and dualization}

Recall that, when $U$, $V$ and $W$ are sets, the \emph{composition}, ${\mathbf
  R}' \ast \mathbf R$, of
relations ${\mathbf R}
\subset U \times V$ and ${\mathbf R'} \subset V \times W$ is
\[ {\mathbf R}' \ast \mathbf R =\{(u, w) \in U \times W\, | \, \exists v \in
V, (u, v) \in \mathbf R \, {\mathrm{and}} \, (v, w)
   \in \mathbf R' \}, \]
and the \emph{inverse} of a relation $\mathbf R \subset U \times V$ is the relation
\[ \overline{\mathbf R} =\{(v, u) \in V \times U\, | \, (u, v) \in \mathbf R\}. \]
We observe that 
\begin{equation}\label{eq01}
\overline{{\mathbf R}' \ast \mathbf R} = \overline{\mathbf R}  \ast  
\overline{{\mathbf R}'}.
\end{equation}
It is clear that, if $\phi : U \to V$ and $\phi' : V \to W$ are maps, and if 
$\mathbf R = \mathrm{\tmop{graph}} \, \phi$ and $\mathbf R' =
\mathrm{\tmop{graph}} \, \phi'$, then
\[ \mathbf R' \ast \mathbf R = \mathrm{\tmop{graph}} (\phi' \circ \phi) . \]
and, if $\phi$ is invertible, 
\[\overline{\mathrm{\tmop{graph}} \, \phi} = \mathrm{\tmop{graph}}
(\phi^{-1}).\]

Let $U$ and $V$ be vector spaces.
The \emph{dual} of a relation $\mathbf R \subset U \times V$
is the relation $\mathbf R^* \subset V^* \times U^*$ defined by
$$
 \mathbf R^* =\{(\beta, \alpha) \in V^* \times U^*\, | \, \langle \alpha, u
\rangle = \langle \beta, v \rangle,
   \forall (u,v) \in \mathbf R \}.
$$
We observe that
\begin{equation}\label{eq02}
\overline{\mathbf R^*} = \overline{\mathbf R}\,^*.
\end{equation}
It is clear that if $\mathbf R = \mathrm{\tmop{graph}} \, \phi$,
where $\phi$ is a linear map
from $U$ to $V$, then $\mathbf R^{*}$ is the graph of the dual
map, $\phi^*$.

When $U$ is a set, we shall call a relation $\mathbf R \subset U
\times U$ a \emph{relation in $U$}.

When $U$ and $V$ are vector bundles over a manifold $M$,
and $\mathbf R \subset U \times V$ is a relation, we denote by
$\underline{\mathbf R}$
the relation on sections induced by $\mathbf R$,
$$
\underline{\mathbf R} = \{ (\underline u, \underline v) \in \Gamma U 
\times \Gamma V \, | \, \forall x \in M,
(\underline u(x),\underline v(x)) \in {\mathbf R} \}.
$$
When $\mathbf R$ is a vector subbundle of $U \times V$, 
then $\underline{\mathbf R}$ is the space of sections of $\mathbf R$. 
If $\mathbf R = \mathrm{\tmop{graph}} \, \phi$, for a vector bundle
morphism $\phi : U \to V$, and if 
${\underline \phi} : \Gamma U 
\to \Gamma V$ is defined by $\forall \underline u \in \Gamma U,
\forall x \in M, {\underline \phi}(\underline u)(x)=\phi(u(x)$,
then $\underline {\mathbf R} = \underline{\mathrm{\tmop{graph}} \, \phi}$.
In the rest of this paper, we shall neglect to underline the notation
for sections and for relations in spaces in sections.

We remark that, when $U$, $V$ and $W$ are vector bundles, and $\mathbf R$ and
$\mathbf R'$ are relations defined by vector subbundles (of constant rank)
of $U \times V$ and $V \times W$ respectively, $\mathbf R' \ast \mathbf R$ is not
a subbundle of $U \times W$ unless a constant rank condition is
satisfied. (See, {\it e.g.}, \cite{W}.) 

\subsection{Nijenhuis relations in Leibniz algebras}
Let $\mathbf N$ be a relation in a Leibniz algebra $(E, [~,~])$. (See
a definition of Leibniz algebras in the Appendix.)
We consider the  real-valued function defined on a subset of 
$E \times E \times E\times E 
\times  E^{\ast} \times  E^{\ast} \times E^*$ by 
\begin{equation} \label{T} 
\begin{array}{lll}
& & \mathbf{T} (\mathbf{N})(u_1, v_1, u_2,
  v_2,\alpha,\alpha', \alpha'')\\ 
&=&\langle \alpha, [v_1, v_2] \rangle -  \langle \alpha', [v_1,
      u_2] + [u_1, v_2] \rangle +  \langle  \alpha'', [u_1, u_2] \rangle,
\end{array}
\end{equation}  
for all $u_1, v_1, u_2, v_2 \in E, \alpha, \alpha', \alpha'' \in
E^*$ such that $(u_1, v_1) \in {{\mathbf N}}, (u_2, v_2) \in {{\mathbf
    N}}$, $(\alpha, \alpha')  \in
  {\mathbf N}^{\ast},( \alpha', \alpha'') \in {\mathbf N}^{\ast}$. 

\begin{definition}
Let $(E, [~,~])$ be a Leibniz algebra, and let $\mathbf N$ be a
relation 
in~$E$. The function $\mathbf{T} (\mathbf{N})$ 
is called the \emph{torsion} of the
relation $\mathbf N$.
A \emph{Nijenhuis relation} in $E$ is a subset $\mathbf{N}$ of $E
\times  E$ 
such that its torsion $\mathbf{T} (\mathbf{N})$ vanishes.
\end{definition}

It is easy to see, using relation \eqref{eq02}, that 
\begin{equation}\label{inverse}
\mathbf{T} (\mathbf{N}) = \mathbf{T} (\overline{\mathbf{N}}).
\end{equation}

We now prove that Nijenhuis relations generalize Nijenhuis tensors.

\begin{proposition}\label{nij}
A linear map $N : E \to E$ is a Nijenhuis tensor
  if and only if $\mathrm{\tmop{graph}} \hspace{0.25em} N$ is a Nijenhuis
  relation in $ E$.
\end{proposition}

\begin{proof}
  The graph of $N$ is the relation, 
$$\mathrm{\tmop{graph}} \hspace{0.25em} N
  =\{(u, N u) \in E \times E\, | \,u \in E\},
$$ 
and its dual is the graph of the dual
$N^{\ast}$ of $N$, 
$$\mathrm{\tmop{graph}} (N^{\ast}) =\{(\alpha, N^{\ast} \alpha)
  \in E^{\ast} \times E^{\ast} \, | \, \alpha \in E^{\ast} \}.$$ 
Thus,
$\mathrm{\tmop{graph}} \hspace{0.25em} N$ is a Nijenhuis relation if and
  only if, for all $u_1, u_2 \in E$, $\alpha \in E^{\ast}$,
  \[ \langle \alpha, [N u_1, N u_2] \rangle 
- \langle N^{\ast} \alpha, [N u_1, u_2] +  [u_1, N u_2] \rangle + \langle
     (N^{\ast})^2 \alpha, [u_1, u_2] \rangle = 0. \]
This is equivalent to $\langle  \alpha, {T}N (u_1, u_2) \rangle = 0$, where
  $TN$ is the Nijenhuis torsion of $N$. By the non-degeneracy of
  the pairing, this condition is equivalent to the vanishing of the Nijenhuis
  torsion of $N$.
\end{proof}

\subsection{The torsion as a relation}
For relations  $\mathbf R \subset U \times V$ and  
${\mathbf R}' \subset V \times W$, we define 
$$ {\mathbf R}' \diamond \mathbf R =\{(u, v,  w) 
\in U \times V \times W\, | (u, v) \in \mathbf R \, \, {\mathrm{and}}
\, \, (v, w) \in \mathbf R'\}.
$$
(Cf. \cite{W}, where ${\mathbf R}' \diamond \mathbf R$ is 
viewed as a fiber product ${\mathbf R'} \times_V \mathbf R$.) Then
the projection of ${\mathbf R}'
\diamond \mathbf R$ on $U \times W$ is ${\mathbf R}' \ast \mathbf R$.

If $\phi : U \to V$ and $\phi' : V \to W$, and if 
$\mathbf R = \mathrm{\tmop{graph}} \, \phi$ and $\mathbf R' =
\mathrm{\tmop{graph}} \, \phi'$, then
$$ {\mathbf R}' \diamond \mathbf R = \{(u, \phi(u), \phi'(\phi(u)) \,
| \, u \in U\},
$$
i.e., $\mathrm{\tmop{graph}} \, \phi' \diamond  \mathrm{\tmop{graph}}
\, \phi
=
 \mathrm{\tmop{graph}} (\phi ,   \phi' \circ \phi) \subset U \times (V
 \times W)$.

When $V=U$, for a relation $\mathbf R \subset U \times U$, 
we set 
$$
\mathbf
R^{(2)} = {\mathbf R}\diamond \mathbf R = 
 \{ (u, u', u'') \in U \times U \times U \, | \, (u,u')\in
\mathbf R \, \, {\mathrm{and}} \, \, (u',u'') \in \mathbf R \}.
$$
With this notation, we see that the vanishing of ${\mathbf{T}
(\mathbf{N})}$ defined by \eqref{T} defines a relation, 
$$\widehat{\mathbf{T}
(\mathbf{N})} \subset (\mathbf{N} \times \mathbf{N}) \times
(\mathbf{N}^{\ast})^{(2)}.
$$

\subsection{The case of Leibniz algebroids}
We consider a Leibniz
algebroid over a manifold $M$ (see the Appendix) and we 
denote the Leibniz (Loday) bracket on sections by $[~,~]$.
Let $\mathbf N$ be a relation in a Leibniz algebroid $(E, \rho, [~,~])$. 
We define the torsion of  $\mathbf N$ by equation 
\eqref{T}, obtaining a function on the subset
$(\mathbf{N} \times \mathbf{N}) \times
(\mathbf{N}^{\ast})^{(2)}$ of $\Gamma E \times \Gamma E \times \Gamma E
\times \Gamma E \times\Gamma (E^*) \times\Gamma (E^*) \times\Gamma
(E^*)$.
 
\begin{proposition}\label{nijoid}
A vector bundle morphism $N : E \to E$ is a Nijenhuis tensor
  if and only if $\mathrm{\tmop{graph}} \hspace{0.25em} N$ defines a Nijenhuis
  relation in $\Gamma E$.
\end{proposition}

\begin{proof}
We introduce the map on sections induced by the 
endomomorphism $N$ of $E$. The proof is then formally analogous to that of
the case of Leibniz algebras.
\end{proof}

\section{Hamiltonian pairs}\label{Hpairs}
Pairs of compatible Poisson structures on Lie algebroids 
were studied in \cite{KSR} using the big
bracket. Here we  study them using relations, following Gelfand 
and Dorfman \cite{GD0} \cite{GD} and Dorfman \cite{D1987} \cite{D}.
In this section and in Sections~\ref{3} and~\ref{4}, 
we shall consider a Lie algebroid $(A,\mu)$ and we shall denote the bracket
of sections of $A$ defined by $\mu$ by $[~,~]$. (See the Appendix for
definitions and notation.)
Nijenhuis relations have just been defined in Section \ref{sec1}.

\subsection{Poisson structures and Hamiltonian pairs}
Recall that $\pi \in \Gamma (\wedge^2 A)$ is a \emph{Poisson structure} (also
called a Hamiltonian structure)
if $[\pi,
\pi] = 0$, where $[~, ~]$ is the Schouten--Nijenhuis bracket of
multivectors defined by $\mu$. 
We denote by the same letter a bivector $\pi$ and 
the map $\pi : A^{\ast} \to A$ defined by
$\pi \xi = i_{\xi} \pi$ for all $\xi\in \Gamma (A^{\ast})$, as well as
the map on sections defined by $\pi$. 
When $\pi \in \Gamma (\wedge^2 A)$, set, for $\xi, \eta
\in \Gamma (A^{\ast})$,
\[ [\xi_1,  \xi_2 ]_{\pi} = L_{\pi \xi_1}  \xi_2  - L_{\pi  \xi_2 }
\xi_1 - d (\pi (\xi_1, \xi_2 )) . \]
The following lemma is well known (see proposition 3.1 of {\cite{KSM}}).
\begin{lemma}\label{lemma2}
For all $\xi_1, \xi_2 \in \Gamma (A^{\ast})$,
\begin{equation}\label{pipi}
\frac12 [\pi,\pi](\xi_1, \xi_2 ) = [\pi \xi_1, \pi\xi_2  ]- \pi[\xi_1,
  \xi_2]_{\pi}.
\end{equation}
A bivector $\pi$ is a Poisson structure if and only if, for
  all $\xi_1, \xi_2\in \Gamma (A^{\ast})$,
\[ [\pi \xi_1, \pi  \xi_2 ] =  \pi [\xi_1,  \xi_2 ]_{\pi} . \]
\end{lemma}

\begin{definition}
  Poisson structures $\pi$ and $\pi'$ on $A$ are said to be 
\emph{compatible} if $\pi + \pi'$ is a Poisson structure.
When Poisson structures $\pi$ and $\pi'$ are compatible,
$(\pi,\pi')$ is said to be a
\emph{bi-Hamiltonian structure} or a \emph{Hamiltonian pair}.
\end{definition}

Poisson structures $\pi$ and $\pi'$ constitute a Hamiltonian
pair if and
only if $[\pi,\pi'] \!  =~\! \! 0$.

\subsection{The relation defined by a Hamiltonian pair}

For bivectors $\pi$ and $\pi'$, set
\begin{equation}\label{eq2}
\mathbf{N}(\pi, \pi') \! = \! \{(x, y) \in A \times A \, | \, \exists \xi
  \in A^{\ast}, (\xi, x) \in \mathrm{\tmop{graph}} \, \pi', (\xi,
  y) \in \mathrm{graph} \, \pi \},
\end{equation}
which is to say
\begin{equation}\label{eq4a} \mathbf{N}(\pi, \pi') =
  \mathrm{\tmop{graph}} \, \pi \ast
   \overline{\mathrm{\tmop{graph}} \, \pi'} .
\end{equation}
Then
\[ \mathbf{N}(\pi, \pi') =\{(\pi' \xi, \pi \xi) \in A \times A\, | \, \xi \in
   A^{\ast} \}, \]
and
\[ \mathbf{N}(\pi, \pi')^{\ast} =\{(\xi, \xi') \in A^{\ast} \times A^{\ast}
   \, | \, \pi \xi = \pi' \xi' \}. \]
\begin{theorem}\label{graphs}
Let $\pi$ and $\pi'$be bivectors. The torsion of the relation
$\mathbf{N}(\pi, \pi')$ satisfies the relation
 \begin{equation}\label{pair}
\begin{array}{lll}
&& 2 \mathbf{T}( \mathbf{N}(\pi, \pi'))(\xi_1,\xi_2, \xi, \xi',\xi'')\\
&=& \langle \xi, [\pi,\pi](\xi_1,\xi_2) \rangle
+\langle \xi'', [\pi',\pi'](\xi_1,\xi_2) \rangle
- 2 \langle  \xi', [\pi,\pi'](\xi_1,\xi_2)\rangle.
\end{array}
\end{equation}
for all $\xi_1, \xi_2, \xi, \xi',\xi'' \in
\Gamma(A^*)$ such that $\pi \xi = \pi' \xi'$ and $\pi \xi' = \pi' \xi''$.
\end{theorem}

\begin{proof}
The proof is based on Lemma \ref{lemma2}.
On the one hand we obtain, using \eqref{pipi}, for all $\xi_1, \xi_2 \in \Gamma(A^*)$,
$$[\pi +\pi', \pi + \pi'](\xi_1,\xi_2) = [\pi,\pi] (\xi_1,\xi_2)
        + [\pi',\pi'](\xi_1,\xi_2)
+ 2 Q(\xi_1,\xi_2),
$$ 
where 
$$ Q(\xi_1,\xi_2)= [\pi \xi_1, \pi' \xi_2] + [\pi'
    \xi_1, \pi \xi_2]  - \pi [\xi_1, \xi_2]_{\pi'} - \pi' [\xi_1,
  \xi_2]_{\pi}
= [\pi,\pi'](\xi_1,\xi_2).
$$
On the other hand, the expression for the torsion of 
$\mathbf{N}(\pi, \pi')$ is, for all $\xi_1, \xi_2, \xi, \xi', \xi'' \in \Gamma
  (A^{\ast})$, such that $\pi \xi = \pi' \xi'$ and $\pi \xi' = \pi' \xi''$,
$$
\begin{array}{lll}
&&\mathbf{T}( \mathbf{N}(\pi, \pi'))(\xi_1,\xi_2, \xi, \xi',\xi'')\\
&=& \langle \xi, [\pi \xi_1, \pi\xi_2] \rangle 
- \langle \xi', ([\pi' \xi_1, \pi\xi_2] + [\pi \xi_1, \pi'\xi_2]) \rangle 
+ \langle \xi'', [\pi' \xi_1, \pi'\xi_2] \rangle . 
\end{array}
$$
Using \eqref{pipi} and the skew-symmetry of $\pi$ and $\pi'$, we
obtain
$$
\begin{array}{lll}
&& 2 \mathbf{T}( \mathbf{N}(\pi, \pi'))(\xi_1,\xi_2, \xi, \xi',\xi'')\\
&=&\langle \xi, [\pi, \pi](\xi_1,\xi_2) \rangle
+
\langle \xi'', [\pi', \pi'](\xi_1,\xi_2) \rangle
- 2 \langle \xi', Q(\xi_1,\xi_2)\rangle.
\end{array}
$$
Whence \eqref{pair}. 
\end{proof}

As an immediate consequence of relation \eqref{pair}, we obtain the
following corollary.

\begin{corollary}\label{2.4}
If $(\pi, \pi')$ is a Hamiltonian pair, then $\mathbf{N}(\pi, \pi')$ is
a Nijenhuis relation.
\end{corollary}

For Poisson structures $\pi$ and $\pi'$, set $K=
\pi^{-1} ({\mathrm {Im}}\,\pi') \cap \pi'^{-1} ({\mathrm {Im}}\,\pi)
\subset~\!A^*$, and let $K^\perp$ be the orthogonal of $K$ in $A$.
Relation \eqref{pair} and the non-degeneracy of the pairing imply the
following proposition.

\begin{proposition}\label{2.5}
(i) If $\pi$ and $\pi'$ are Poisson structures such that 
$\mathbf{N}(\pi, \pi')$ is a Nijenhuis relation,
then, for all $\xi_1,\xi_2 \in \Gamma(A^*)$,
$[\pi,\pi'](\xi_1,\xi_2)$ is an element of $K^\perp$. 

(ii) If, in addition, $K = A^*$, then
$(\pi, \pi')$ is a Hamiltonian pair.
\end{proposition}

We introduce the following convenient definition.

\begin{definition}\label{Ppair}
Poisson bivectors $\pi$ and $\pi'$ on $A$ are said to be 
a \emph{Poisson pair} if $\mathbf{N}(\pi, \pi')$ is a Nijenhuis relation.
A \emph{non-degenerate Poisson pair} is a Poisson pair in which
both bivectors are non-degenerate. 
\end{definition}

If $\pi$ and $\pi'$ are both
non-degenerate bivectors, $K = A^*$. Therefore we obtain the following
corollary of Proposition \ref{2.5}.

\begin{corollary}\label{Ppairbis} 
Any Hamiltonian pair is a
Poisson pair, and, conversely, if $(\pi,\pi')$
is a non-degenerate Poisson pair, 
then $(\pi, \pi')$ is a Hamiltonian pair.
\end{corollary}

It follows from the definition (formula \eqref{eq4a}) and equation \eqref{eq01} that 
\begin{equation}\label{symm} 
\mathbf{N}(\pi', \pi) =
\overline{ \mathbf{N}(\pi, \pi')}.
\end{equation}
From relation
\eqref{inverse}, we find that the torsion is
symmetric in $\pi$ and $\pi'$, 
\begin{equation}
\mathbf{T}(\mathbf{N}(\pi',\pi)) = \mathbf{T}(\mathbf{N}(\pi,\pi')).
\end{equation}
(This fact also follows from formula
\eqref{pair} of Theorem \ref{graphs}, since $[\pi, \pi'] = [\pi',
  \pi]$).
Therefore $(\pi, \pi')$ is a Poisson pair
if and only if  $(\pi', \pi)$ is a Poisson pair. 

\medskip

When $\pi$ is non-degenerate, 
\[ \mathbf{N}(\pi', \pi) = \mathrm{\tmop{graph}} \pi' \ast
  {\mathrm{\tmop{graph}} (\pi^{-1})} = \mathrm{\tmop{graph}}(\pi'
   \circ \pi^{-1}).\]
Corollaries \ref{2.4} and \ref{Ppairbis} and Proposition \ref{nij} 
imply the following well known result. (See, {\it e.g.}, 
\cite{KSR}.)

\begin{proposition}
(i) Assume that $(\pi,\pi')$ is a Hamiltonian pair, where $\pi$ is
  non-degenerate. Then $N =
  \pi' \pi^{- 1}$ is a Nijenhuis tensor.

(ii) Assume that $\pi$ and $\pi'$ are non-degenerate Poisson structures and
that $N =
  \pi' \pi^{- 1}$ is a Nijenhuis tensor. Then $(\pi, \pi')$ is a
  Hamiltonian pair.
\end{proposition}

\begin{remark}\label{remark}\rm{
  For any Lie algebroid $(A,\mu)$, 
$(\Omega^\bullet = \Gamma (\wedge^\bullet A^{\ast}), \d)$ is a
  complex over $\Gamma A$ in the sense of Dorfman (\cite{D}, p. 11), 
with a non-degenerate
  pairing of $\Gamma A$ and $\Omega^1 = \Gamma (A^{\ast})$. Therefore
  Proposition \ref{nijoid} above
is a particular case of proposition
  3.15 of Dorfman {\cite{D}}. It is claimed in 
theorem 3.16 of \cite{D} that, when $\pi$ and $\pi'$ are Poisson
structures, the vanishing of the torsion of the
relation ${\mathbf{N}}(\pi,\pi')$ implies that $(\pi, \pi')$ is a
Hamiltonian pair,
but the proof relies on an application of its
proposition 3.11 that neglects the condition (in our notation) $\xi'
\in K$.}
\end{remark}

\begin{remark}\rm{
  \label{th2}Let $\pi$ and $\pi'$ be Poisson bivectors on $A$. Then $(\pi,
  \pi')$ is a Hamiltonian pair if and only if $(\lambda \pi,
  \lambda' \pi')$ is a Hamiltonian pair for all $\lambda \in \mathbb{R}$ and all $\lambda' \in
  \mathbb{R}$.
This result follows from the bilinearity of the Schouten--Nijenhuis bracket.
We can also show that the assumption that $(\pi, \pi')$ is a
Poisson pair implies that
$\mathbf{N}(\lambda \pi, \lambda' \pi')$ is a Nijenhuis relation.
Thus $(\pi, \pi')$ is a Poisson pair if and only if $(\lambda \pi,
  \lambda' \pi')$ is a Poisson pair 
for all $\lambda, \lambda' \in \mathbb{R}$.}
\end{remark}

\section{Dirac pairs, Poisson pairs and 
presymplectic pairs}\label{3}

\subsection{Dirac pairs}

Let $A$ be 
a vector bundle, and let $A^*$ be the dual vector bundle. For
relations $L \subset A \times A^*$ and $L' \subset A \times A^{\ast}$,
we consider the relation in $A$,
\begin{equation}\label{reldirac}
 \mathbf{N}_{L, L'} =\{(x, y) \in A \times A\, | \, \exists \xi \in A^{\ast}, (x,
   \xi) \in L', (y, \xi) \in L\}.
\end{equation} 
It is clear that
\begin{equation}\label{defpair}
\mathbf{N}_{L, L'} = \overline L * L'.
\end{equation} 
We observe that 
$$
\mathbf{N}_{L', L} =\overline{\mathbf{N}_{L, L'}}.
$$

Assume that $(A,\mu)$ 
is a Lie algebroid, and that $E
  = A \oplus A^{\ast}$ is equipped with the Dorfman bracket. (See the
  Appendix.)
Recall that a 
maximally isotropic subbundle of $ A \oplus A^{\ast}$ whose space of
sections is closed under the Dorfman bracket is called a \emph{Dirac
structure} on $A$.

Following Dorfman \cite{D1987} \cite{D}, we introduce the
following compatibility condition on Dirac structures.

\begin{definition}\label{Dpair} 
 Dirac structures
  $L$ and $L'$ on $A$ are said to be a \emph{Dirac pair}
  if $\mathbf{N}_{L, L'}$ defined by \eqref{reldirac} 
is a Nijenhuis relation in $A$.
\end{definition}

This definition can be generalized: relations $L$ and $L'$ 
in $A \oplus A^*$ are said to be 
compatible if $\mathbf{N}_{L, L'} = \overline L * L'$ is a Nijenhuis
relation in $A \oplus A^*$.

\subsection{Poisson pairs}
A bivector on $A$ defines a Poisson structure if
and only if its graph is a Dirac structure on $A$. (See \cite{LWX}.)

It follows from definitions \eqref{eq2} and \eqref{reldirac} that, 
when $\pi$ and $\pi'$ are bivectors, if $L =
\overline
{\mathrm{\tmop{graph}} \, \pi}$ and $L' = \overline{\mathrm{\tmop{graph}} \, \pi'}$,
then 
$$
\mathbf{N}_{L, L'} = \mathrm{\tmop{graph}} \,
  \pi * \overline{\mathrm{\tmop{graph}} \, \pi'} =\mathbf{N}(\pi, \pi') .
$$ 
\begin{theorem}
(i) Bivectors $\pi$ and $\pi'$ constitute a Poisson pair if and only
if their graphs constitute a Dirac pair.

\noindent (ii) If $(\pi,\pi')$ is a Hamiltonian pair,
then $(\overline{\mathrm{\tmop{graph}} \, \pi},\overline{\mathrm{\tmop{graph}} \, \pi'})$
is a Dirac pair.

\noindent (iii) 
Conversely, if
$(\overline{\mathrm{\tmop{graph}} \, \pi},\overline{\mathrm{\tmop{graph}} \, \pi'})$ 
is a Dirac pair
and if  $\pi$ and $\pi'$ are non-degenerate bivectors, then
$(\pi,\pi')$ is a Hamiltonian pair.
\end{theorem}

\begin{proof}
(i) follows from Definitions \ref{Ppair} and \ref{Dpair},
(ii) is a consequnce of Corollary \ref{2.4}, and (iii) is a
  consequence of Proposition \ref{2.5}.
\end{proof}

\subsection{Presymplectic pairs}
Dually, we can define presymplectic pairs.
Recall that a \emph{presymplectic structure} on $A$ is defined by
a closed $2$-form on
$A$. We denote by the same letter a $2$-form
$\omega \in\Gamma(\wedge^2A^*)$ and the map
$\omega : A \to A^*$ defined by $\omega x = - i_x\omega$ for all $x \in A$,
as well as the map on sections induced by $\omega$.

A $2$-form on $A$ defines a presymplectic structure if
and only if its graph is a Dirac structure on $A$. (See \cite{LWX}.)
\begin{definition}
  If $\omega$ and $\omega'$ are presymplectic structures whose graphs
constitute a
  Dirac pair, $(\omega, \omega')$ is called a \emph{presymplectic pair}.
If, in addition, $\omega$ and $\omega'$ are non-degenerate, $(\omega, \omega')$ is called a \emph{symplectic pair}.
\end{definition}

For $L= \mathrm{\tmop{graph}} \, \omega$, 
$L'= \mathrm{\tmop{graph}} \, \omega'$,
$$
\mathbf{N}_{L,L'}= \overline{\mathrm{\tmop{graph}} \, \omega} *
\mathrm{\tmop{graph}} \, \omega'.
$$

Presymplectic pairs were introduced by Dorfman \cite{D} under the name
`symplectic pairs'. (She called `symplectic' what we call
`presymplectic', and `invertible symplectic' what we call
`symplectic'.)

In Section \ref{Hpairs}, Corollary \ref{Ppairbis},
we observed 
that the non-degenerate Poisson pairs coincide with
the non-degenerate Hamiltonian pairs. 
The following result relates
non-degenerate Poisson pairs to symplectic pairs.

\begin{theorem}\label{1-1}Symplectic pairs are in one-to-one correspondence 
with non-degenerate Poisson pairs. 
\end{theorem}

\begin{proof}
In fact, if $\pi \! =  \omega^{-1}$ and $\pi' =  \omega'^{-1}$,
$\overline{\mathrm{\tmop{graph}} \, \omega} \, * \, 
\mathrm{\tmop{graph}} \, \omega' = \mathrm{\tmop{graph}} \,
(\omega^{-1} \circ~\omega')\break
= \mathrm{\tmop{graph}} \, (\pi \circ \pi'^{-1}) =
\mathrm{\tmop{graph}} \, \pi \ast
   \overline{\mathrm{\tmop{graph}} \, \pi'} = {\mathbf{N}}(\pi,\pi')$.
Therefore $(\omega,\omega')$ is a symplectic pair if and only if
$(\pi,\pi')$ is a Poisson pair.
\end{proof}

\begin{example}\label{example}\rm{
Examples of presymplectic pairs arise in the theory of Monge-Amp\`ere
operators. (See \cite{LRC}, \cite{KLR}, chapters 6 and 20, 
and \cite{KSR}, section 13.) Let $M = T^* \mathbb R^2$
and let $\Omega$ be the canonical symplectic form on
$M$. Here $A = TM$. 
In  canonical coordinates $(q^1,q^2,p_1,p_2)$ on $M$,
 $\Omega = \d q^1\wedge \d p_1+\d q^2\wedge \d p_2$. 
Explicit examples of presymplectic pairs 
$(\Omega, \omega)$ are defined by 
$$\omega = \omega_H= \d q^1\wedge \d p_1-\d q^2\wedge \d p_2,$$
$$\omega = \omega_E = \d q^1\wedge \d p_2-\d q^2\wedge \d p_1,$$
$$\omega = \omega_P = \d q^1\wedge \d p_2.$$
The $2$-form $\omega_H$ (resp., $\omega_E$) is a closed, normalized,
effective $2$-form on $M$ in the sense of \cite{LRC}
\cite{KLR}, corresponding to a hyperbolic (resp., elliptic)
Monge-Amp\`ere
equation with constant coefficients. In fact, for $\omega= \omega_H$
or $\omega= \omega_E$, $\omega \wedge \Omega = 0$, {\it i.e.}, $\omega$ is
effective, while 
 $\omega \wedge \omega = - \Omega \wedge \Omega$ if $\omega=
\omega_H$ , and $\omega
\wedge \omega = \Omega \wedge \Omega$ if $\omega = \omega_E$, 
{\it i.e.}, $\omega$ is normalized. 
The non-degenerate $2$-forms $\omega_H$ and $\omega_E$ give rise to symplectic
pairs, but $\omega_P$, corresponding to the parabolic case characterized by
$\omega \wedge \omega = 0$, gives rise to a presymplectic pair. 
It is easy to prove these facts in the present framework.
Let us denote the bivector inverse of $\Omega$ (resp., $\omega_H$,
$\omega_E$)
by $\pi_\Omega$ (resp., $\pi_H$, $\pi_E$). 
Since each of these bivectors has constant
coefficients, $(\pi_\Omega, \pi_H)$ and $(\pi_\Omega, \pi_E)$ constitute
non-degenerate Poisson pairs, and this fact implies that $(\Omega,
\omega_H)$ and $(\Omega,
\omega_E)$ are symplectic pairs. The fact that $(\Omega,
\omega_P)$ is a presymplectic pair cannot be proved by a similar
argument since $\omega_P$ is not invertible.  
We calculate the torsion of the
$(1,1)$-tensor 
$\pi_\Omega \circ \omega_P$ and we see that it vanishes.} 

If $\omega$ is a closed, effective $2$-form with positive Pfaffian, then
$(\Omega, \omega)$ 
is a `symplectic couple' in the
sense of \cite{G}. 
If, in addition, $\omega$ is
normalized so that its Pfaffian is equal to~$1$ -- which is the case for
$\omega = \omega_E$ --
the pair $(\Omega, \omega)$ is then a `conformal symplectic couple' as
defined in \cite{G}. 
Any conformal symplectic couple on a $4$-manifold $M$
defines a Dirac pair on $TM$.
See \cite{G} for criteria for the existence of symplectic couples
and of conformal symplectic couples on $4$-manifolds. 
\end{example}

\section{P$\Omega$- and $\Omega$N-structures}\label{4} 
We now characterize P$\Omega$- and $\Omega$N-structures in terms of Dirac pairs.

\subsection{P$\Omega$-structures} 
\begin{definition}
A bivector $\pi$ and a $2$-form $\omega$ define a P$\Omega$-structure 
on a Lie algebroid $(A,\mu)$
if $\pi$ is a Poisson bivector, and both $\omega$ and $\omega_N$ are 
closed,
where $N = \pi \circ \omega$ and $\omega_N = \omega \circ N$.
\end{definition}

We consider the relation $\mathbf{N}_{L, L'}$ in
the case of 
$L = \overline{\mathrm{\tmop{graph}} \hspace{0.25em} \phi}$, $L'~=~\mathrm{\tmop{graph}}
\hspace{0.25em} \phi'$, where $\phi : A^{\ast} \to A$ and $\phi' : A \to
A^{\ast}$. By definition \eqref{reldirac},
$$     \mathbf{N}_{L, L'} 
    =  \{(x, y) \in A \times A\, | \, \exists \xi \in A^{\ast},
     \xi = \phi' x, y = \phi \xi\},
$$
and therefore
\begin{equation}
  \label{comp} \mathbf{N}_{L, L'} =\{(x, y) \in A \times A\, | \,y = \phi (\phi'
  x)\}= \mathrm{\tmop{graph}} \hspace{0.25em} (\phi \circ \phi') .
\end{equation}

\begin{proposition}Let $\pi$ be a Poisson bivector and let $\omega$ be a
closed $2$-form.
Then 
$(\overline{\mathrm{\tmop{graph}} \hspace{0.25em} \pi},
\mathrm{\tmop{graph}} \hspace{0.25em} \omega)$ is a Dirac pair if and only if
$\pi \circ \omega$ is a Nijenhuis tensor.
\end{proposition}

\begin{proof}
We apply formula
\eqref{comp} for $L = \overline{\mathrm{\tmop{graph}} \hspace{0.25em}
\pi}$ and $L'= \mathrm{\tmop{graph}} \hspace{0.25em} \omega$. The
result then follows from Proposition \ref{nij}.
\end{proof}

If $N$ is a $(1,1)$-tensor on a Lie algebroid $(A,\mu)$, 
we define the operator on forms $i_N = \{N, \cdot\}$, where $\{~,~\}$
is the big bracket (see the Appendix),
and we let $\d_N$ be the graded commutator, $
\d_N = [i_N, \d]$.
Then, if $\omega$ is closed, $$\d\omega_N=\d_N\omega.$$
Therefore 
a Poisson bivector $\pi$ and a closed $2$-form $\omega$ on $(A,\mu)$
define a $P\Omega$-structure if and only if $\d_N \omega = 0$, where
$N = \pi \circ \omega$.
The proof of the following theorem relies on formula (7.4) in
\cite{KSR}.

\begin{theorem}\label{oldtheorem}
(i)  If a Poisson structure $\pi$ and a presymplectic
  structure $\omega$ 
  constitute a $P \Omega$-structure, their graphs constitute a Dirac
  pair.
  
\noindent (ii) Conversely, if the graphs of a Poisson structure $\pi$ and a presymplectic
  structure $\omega$ constitute a Dirac pair, and if $\pi$ is non-degenerate,
  then $\pi$ and $\omega$ constitute a $P \Omega$-structure.
\end{theorem}

\begin{proof}
  Let $\pi$ (resp., $\omega$) be a Poisson (resp., presymplectic) structure.
  If $(\pi, \omega)$ constitutes a $P \Omega$-structure, then it follows from
  Corollary 4.4 of {\cite{A}} or Theorem 8 of {\cite{KSR}} that $N = \pi \circ
  \omega$ is a Nijenhuis tensor. Therefore the graph of $N$ is a Nijenhuis relation,
  and $(\pi, \omega)$ constitutes a Dirac pair.
  
  Conversely, we must prove that if $\pi$ is a non-degenerate
  Poisson structure and $\omega$ is a closed $2$-form, 
and $\pi \circ  \omega$ is a Nijenhuis tensor,
then $\d_N \omega =~\!0$. Firstly, it follows from formula
  (7.4) of {\cite{KSR}} that, when $\pi$ is a Poisson bivector,
$\omega$ is a
  closed $2$-form and $N = \pi \circ \omega$ is a Nijenhuis tensor on
  $A$,  then $\{\pi, \d_N \omega\}= 0$. 
Secondly, we must prove that if $\pi$ is
  invertible, then $\{\pi, \d_N \omega\}= 0$ implies that $\d_N \omega =
  0$. For any bivector~$\pi$,
  $2$-form $\sigma$ and $k$-form $\alpha$, $k$ a nonnegative integer, 
 if $\sigma$ is the
  inverse of $\pi$, 
applying the Jacobi identity to $\{\sigma, \{\pi, \alpha\}\}$ yields
  \[ \{\sigma, \{\pi, \alpha\}\}=\{\{\sigma, \pi\}, \alpha\}\}+\{\pi,
     \{\sigma,  \alpha\}\}=\{\{\sigma, \pi\}, \alpha\}\},\]
since the big bracket of any two forms vanishes. 
When $\sigma$ is the inverse of $\pi$, then $\{\sigma, \pi\} = -
\mathrm{Id}_A$. By (2.4) of \cite{KSR}, for any $k$-form $\alpha$, 
$\{\mathrm{Id}_A,  \alpha\}= k \alpha$. We conclude that, for a
non-degenerate bivector $\pi$, $\{\pi, \d_N
\omega\}= 0$ implies that $\d_N \omega = 0$,
which proves (ii).
\end{proof}

\subsection{$\Omega$N-structures}

Let $N$ be a $(1,1)$-tensor and $\omega$ a $2$-form on $(A, \mu)$
such that 
$ \omega  \circ N  = N^* \circ \omega$.
Then $\omega_N$ defined by $\omega_N = \omega \circ N$ is 
a $2$-form.

\begin{definition}
A $2$-form $\omega$ and a $(1,1)$-tensor $N$ define an
\emph{$\Omega$N-structure} on a Lie algebroid $(A,\mu)$
if $\omega \circ N  = N^* \circ \omega$,
$N$ is a Nijenhuis tensor, and both $\omega$ and $\omega_N$ are 
closed,
where $\omega_N = \omega \circ N$.
\end{definition}

\begin{example}\label{exampleomegaN}
\rm{In the notation of Example \ref{example}, in 
  coordinates $(q^1,q^2,p_1,p_2)$, let $N_H= \Omega^{-1} \circ
  \omega_H$ and 
$N_E= \Omega^{-1} \circ \omega_E$, so that 
$$N_H= \begin{pmatrix}
1&0&0&0\\
0&-1&0&0\\
0&0&1&0\\
0&0&0&-1
\end{pmatrix} \quad {\mathrm{and}} \quad N_E= \begin{pmatrix}
0&-1&0&0\\
1&0&0&0\\
0&0&0&1\\
0&0&-1&0
\end{pmatrix}.
$$
Then $(\Omega, N_H)$ and $(\Omega, N_E)$ are
$\Omega N$-structures on $T(T^*{\mathbb R}^2)$, with $N_H^2= {\mathrm
    {Id}}$ and $N_E^2= - {\mathrm
    {Id}}$. Thus $N_E$ is a complex structure, and $N_H$ is a
product structure on $T^*({\mathbb{R}}^2)$.

Let $N_P= \Omega^{-1} \circ \omega_P$, so that $N_P= \begin{pmatrix}
0&0&0&0\\
1&0&0&0\\
0&0&0&1\\
0&0&0&0\end{pmatrix}$. Then $(\Omega, N_P)$ is an 
$\Omega N$-structure with $N_P^2=0$, so that $N_P$ is a tangent structure.}
\end{example}

\smallskip

We first consider the case where the $2$-form is non-degenerate.

\begin{proposition}\label{prop4.6} 
Let $\omega$ be a non-degenerate $2$-form and $N$ a $(1, 1)$-tensor such that
$\omega_N = \omega \circ N$ is skew-symmetric. 
Then $(\omega, N)$ is an $\Omega N$-structure if and only if
$( \mathrm{\tmop{graph}} \hspace{0.25em} \omega, \mathrm{\tmop{graph}}
  \hspace{0.25em} \omega_N)$ is a Dirac pair.
\end{proposition}

\begin{proof}
When $L= \mathrm{\tmop{graph}} \hspace{0.25em} \omega$ and 
$L'= \mathrm{\tmop{graph}} \hspace{0.25em} \omega_N$,
\[ \mathbf{N}_{LL'} =\{(x,y) \in A \times A
   \, | \, \omega_N x = \omega y \}. \]
Therefore, when $\omega$ is invertible, 
\[ \mathbf{N}_{LL'}=\mathrm{\tmop{graph}} \hspace{0.25em} N .\]   
If $(\omega,N)$ is an $\Omega$N-structure, the graphs of $\omega$ and
$\omega_N$ are Dirac structures and $\mathbf{N}_{LL'}$ is a Nijenhuis
relation, so $( \mathrm{\tmop{graph}} \hspace{0.25em} \omega,
\mathrm{\tmop{graph}}
\hspace{0.25em} \omega_N)$ is a Dirac pair.
Conversely, if $( \mathrm{\tmop{graph}} \hspace{0.25em} \omega,
\mathrm{\tmop{graph}}
\hspace{0.25em} \omega_N)$ is a Dirac pair, $\omega$ and $\omega_N$
 are both closed and the $(1,1)$-tensor $N$ is a Nijenhuis tensor.
\end{proof}

\begin{example}
{\rm The pairs  $( \mathrm{\tmop{graph}} \hspace{0.25em}
  \Omega  \hspace{0.25em} ,  \hspace{0.25em}
\mathrm{\tmop{graph}}
\hspace{0.25em} \omega_H)$, $( \mathrm{\tmop{graph}}  \hspace{0.25em}
\Omega \hspace{0.25em} , \hspace{0.25em}
\mathrm{\tmop{graph}}
\hspace{0.25em} \omega_E)$ and $( \mathrm{\tmop{graph}} \hspace{0.25em} \Omega,
\mathrm{\tmop{graph}}
\hspace{0.25em} \omega_P)$
 are the Dirac pairs associated with the $\Omega$N-structures
 described in Example \ref{exampleomegaN}.}
\end{example}

In the next theorem, the $2$-form $\omega$ is not assumed to be
non-degenerate. A closely related result was proved in \cite{D}, p. 54.
Let $\omega$ be a $2$-form and $N$ a $(1, 1)$-tensor such that
$\omega_N = \omega \circ N$ is skew-symmetric. 

We shall call $(\omega, N)$ a {\it weak $\Omega$N-structure} if $\omega$ and
$\omega_N$ are closed $2$-forms, 
and the torsion of $N$ takes values in the kernel
of $\omega$, {\it i.e.}, 
$\omega(TN(x_1,x_2))=0$, for all sections $x_1, x_2$ of A.

We set ${\mathbf N} = {\mathbf N}_{LL'} = \{(x,y) \in A
\times A \, | \, \omega_N x = \omega y \}$ and 
$$ 
{\mathbf N}^+=
\{ (\omega x, \omega_N x) \in A^* \times A^* \, | \, x \in A \}.
$$
The relation $ {\mathbf N}^+$ is the restriction of the graph of $N^*$
to the image of $\omega$, and a
subset of $ {\mathbf N}^*$. 

\begin{theorem}\label{newtheorem} 
(i) If $(\omega, N)$ is an $\Omega N$-structure, and if $ {\mathbf N}^+
={\mathbf N}^*$, then
$( \mathrm{\tmop{graph}} \hspace{0.25em} \omega, \mathrm{\tmop{graph}}
  \hspace{0.25em} \omega_N)$ is a Dirac pair.

\noindent (ii) If $( \mathrm{\tmop{graph}} \hspace{0.25em} \omega,
\mathrm{\tmop{graph}} \hspace{0.25em} \omega_N)$ is a Dirac pair, then
$(\omega, N)$ is a weak $\Omega N$-structure.
\end{theorem}

\begin{proof}
(i) We must prove that, under the hypotheses of part (i) of 
the theorem, ${\mathbf N}$ is a
Nijenhuis relation. Let $\omega' = \omega_N$ and $\omega'' =
\omega'_N= \omega_{N^2}$. When we evaluate $\d\omega$,
$\d\omega'$ and $\d\omega''$ on triples of vectors as indicated below, 
most terms
cancel out, and we obtain, for all $(x_1,y_1) \! \in \!{\mathbf N}$, $(x_2,y_2) 
\in {\mathbf N}$, and for all $x \in \Gamma A$,
\begin{equation}\label{omegaN}
\begin{array}{lll}
\d \omega(y_1,y_2, x) - \d \omega'(y_1,x_2, x) - \d \omega'(x_1,y_2, x)
+ \d \omega''(x_1,x_2, x)\\
= \langle \omega x ,  [y_1,y_2] \rangle
- \langle \omega' x , [x_2,y_1] + [x_1,y_2] \rangle
+ \langle \omega'' x , [x_1,x_2] \rangle.
\end{array}
\end{equation}
If, in particular, we let $y_i = N x_i$, for $i=1,2$, formula \eqref{omegaN} becomes 
\begin{equation}\label{omegaN1}
\begin{array}{lll}
\d \omega(y_1,y_2, x) - \d \omega'(y_1,x_2, x) - \d \omega'(x_1,y_2, x)
+ \d \omega''(x_1,x_2, x)\\
= \langle \omega x ,  TN(x_1,x_2) \rangle.
\end{array}
\end{equation}
Equation \eqref{omegaN1}
shows that if $(\omega,N)$ is an $\Omega$N-structure, then
$\d\omega''=0$. Thus, if $(\omega, N)$ is an $\Omega N$-structure,
the three $2$-forms, $\omega$, $\omega'$ and
$\omega''$, are closed. 
Therefore, by equation \eqref{omegaN}, if $(\omega,N)$ is an $\Omega$N-structure,
 for all $(x_1,y_1) \in {\mathbf N}$, $(x_2,y_2) \in {\mathbf N}$,
 and for all $x \in \Gamma A$,
\begin{equation}\label{omegaN'}
\langle \omega x ,  [y_1,y_2] \rangle
- \langle \omega' x , [x_2,y_1] + [x_1,y_2] \rangle
+ \langle \omega'' x , [x_1,x_2] \rangle = 0.
\end{equation}
If we assume that ${\mathbf N}^+ = {\mathbf N}^*$, equation
\eqref{omegaN'} expresses the fact that
the torsion of relation $ {\mathbf N}$ vanishes.

(ii) If $( \mathrm{\tmop{graph}} \hspace{0.25em} \omega,
\mathrm{\tmop{graph}} \hspace{0.25em} \omega_N)$ is a Dirac pair,
then both $\omega$ and $\omega_N$ are closed $2$-forms. 
Since ${\mathbf N}^+ \subset {\mathbf N}^*$, the vanishing of the
torsion of ${\mathbf N}$ implies that \eqref{omegaN'} is satisfied 
for all $(x_1,y_1) \in {\mathbf N}$, $(x_2,y_2) \in {\mathbf N}$,
 and for all $x \in \Gamma A$.
If, in particular, $y_1 = N x_1$ and $y_2 = N x_2$, we obtain, for all
$x \in \Gamma A$,
\begin{equation}\label{eqomega'}
\langle \omega x , TN(x_1,x_2)\rangle = 0.
\end{equation}
By the skew-symmetry of $\omega$, formula \eqref{eqomega'} proves (ii).
\end{proof}

Clearly, when $\omega$ is non-degenerate, ${\mathbf N}^+ = 
\mathrm{\tmop{graph}} (N^*) = {\mathbf N}^*$ and we
recover Proposition \ref{prop4.6}.

\begin{remark}{\rm
Equation \eqref{omegaN1} can be generalized to prove that, 
given a closed $2$-form $\omega$ and a Nijenhuis tensor $N$, if
$\omega \circ N$ is closed, then all $2$-forms  $\omega \circ N^2$,
$\omega \circ N^3$, \ldots, $\omega \circ N^p, \ldots$ are closed. 
This property is the basis of the construction of 
a sequence of integrals in involution
for bi-Hamiltonian systems, and for the extension of this property to
systems associated to a Dirac pair \cite{D} \cite{BDSK}}.
\end{remark} 
\bigskip

In conclusion, Dirac pairs constitute a very
general framewok in which Dorfman was able to formulate a
generalization of the Magri--Lenart scheme and therefore important
applications to the study of integrable systems.
Poisson pairs, 
presymplectic pairs, all P$\Omega$-structures, and those
$\Omega$N-structures in which the additional condition ${\mathbf N}^+=
{\mathbf N}^*$ is satisfied,
including the symplectic-Nijenhuis structures, furnish
examples of Dirac structures. 
In~{\cite{BDSK}}, Barakat, De Sole and Kac, working in the framework of
the formal calculus of variations, have constructed an example
of a Dirac pair on an infinite-dimensional functional space
which they used to prove the complete integrability of the non-linear
Schr\"odinger hierarchy. 
Dorfman herself \cite{D} applied her theory to many instances of
integrable systems.

\section{Appendix}\label{appendix}

\subsection{Lie algebroids and proto-bialgebroids}\label{A.1}

A \emph{Leibniz algebra} (also called \emph{Loday algebra})
is a vector space, $E$,  
over a field of characteristic $0$ 
with a bilinear bracket, $[~,~]$, satisfying the Jacobi
identity, $[u,[v,w]]= [[u,v],w] + [v, [u,w]]$ for all $u,v,w$ in $E$.
(The bracket is not assumed to be skew-symmetric.)

When $A \to M$ is a vector bundle, 
let $A[n]$ be the graded manifold obtained from $A$ by
assigning 
degree $0$ to the coordinates on the base and degree $n$ ($n$ a
nonnegative integer) to
the coordinates on the fibers.
Let $\mathcal F$ be the bigraded 
commutative algebra of smooth functions on $T^*[2]A[1]$.
Local coordinates 
on $T^*[2]A[1]$, and their bidegrees are:
$$
\begin{matrix}
x^i& \xi^a& p_i& \theta_a\\
(0,0)& (0,1) & (1,1)&
(1,0)
\end{matrix} 
$$
As the cotangent bundle of a graded manifold, $T^*[2]A[1]$ is
canonically equipped with an even Poisson structure \cite{R2002}. 
Denote by $\{~,~\}$ the even Poisson bracket on $\mathcal F$, which we
call the \emph{big bracket}. (See \cite{yks1992} 
for the case where $M$ is a point, 
$A$ is just a vector space, and 
then $\mathcal F= \wedge^\bullet(A \oplus A^*)$.)

The big bracket is of bidegree
$(-1,-1)$, is skew-symmetric,
$$\{u,v\} = - (-1)^{|u|\,|v|} \{v,u\},
$$
and satisfies the Jacobi
identity,
$$\{u,\{v,w\}\}= \{\{u,v\},w\} + (-1)^{|u| \, |v|} \{v, \{u,w\}\}  ,
$$
for all $u$, $v$, $w \in {\mathcal F}$.
In local coordinates,
$\{ x^i,p_j\} = \delta^i_j$ and $\{\xi^a,\theta_b \} = \delta^a_b$.

A \emph{Lie algebroid} structure 
on $A \to M$ is an element $\mu$ of
$\mathcal F$ of bidegree $(1,2)$ such that
$$
\{\mu,\mu\} = 0  .
$$

By a vector (resp., multivector) on $A$ we mean 
a section of $A$ (resp., $\wedge^\bullet A$).
The Schouten--Nijenhuis bracket of
multivectors $X$ and $Y$ on $A$ is
$$
[X,Y] =  \{\{X,\mu\},Y\}.
$$
In particular, the Lie bracket of $X$, $Y \in {\rm{\Gamma}} A$, is 
$$
[X,Y] = \{\{X,\mu\},Y\} 
$$
and 
the anchor of $A$, $\rho: A \to TM$, is defined by
$$
\rho(X)f = [X,f] =\{\{X,\mu\},f\} ,
$$
for all $X \in {\rm{\Gamma}} A$, $f \in C^{\infty}(M)$.

The operator $\d = \{ \mu, . \}$  is a differential on 
$\Gamma(\wedge^\bullet A^*)$ 
which defines the 
Lie algebroid cohomology of $A$, generalizing 
both the Chevalley--Eilenberg cohomology (when $M$
is a point, $A$ is a Lie algebra), 
and the de Rham cohomology (when $A = TM$).
By a form on $A$, we mean a section of $\wedge^\bullet A^*$. By a
closed form we mean a section $\alpha$ of $\wedge^\bullet A^*$ such
that $\d\alpha =0$.

In a \emph{Leibniz algebroid} the space of sections is a Leibniz
algebra.

A \emph{Lie bialgebroid} 
is defined by elements of $\mathcal F$, $\mu$ of bidegree $(1,2)$
and $\gamma$ of bidegree $(2,1)$, such that $\{\mu +
\gamma , \mu + \gamma \} = 0$.
 More generally, 
a \emph{proto-bialgebroid} 
is defined by elements of $\mathcal F$, $\mu$ of bidegree $(1,2)$, 
$\gamma$ of bidegree $(2,1)$, $\phi \in \Gamma(\wedge^3A)$ of
bidegree $(3,0)$, and $\psi \in \Gamma(\wedge^3A^*)$ of
bidegree $(0,3)$, such 
that $$\{\phi +
\gamma + \mu +\psi, \phi + \gamma + \mu + \psi\} = 0.
$$

Let $(A, \mu, \gamma, \phi, \psi )$ be a proto-bialgebroid. Consider the
\emph{Dorfman bracket} on the sections of
$A \oplus A^*$
defined by
$$
[u,v] =  \{\{u, \phi + \gamma + \mu + \psi \},v\}  ,
$$
for all sections $u$ and $v$ of $A \oplus A^*$. 
Then, $\Gamma(A \oplus A^*)$, 
equipped with the Dorfman bracket, is a Leibniz
algebra and $A \oplus A^*$ is a Leibniz algebroid. 
The skew-symmetrized Dorfman bracket is called
the {\it Courant bracket}, and $A \oplus A^*$ is a {\it Courant
  algebroid}, called the  {\it double}
of the proto-bialgebroid $(A, \mu, \gamma, \phi, \psi )$.

In particular, if $(A,\mu)$ is a Lie algebroid, 
considered as a trivial Lie bialgebroid,
the Dorfman bracket on $A \oplus A^*$ is explicitly, 
$$[X+\xi, Y +\eta] = [X,Y] + {\mathcal
      L}_X\eta-i_Y(\d\xi),
$$
for all vectors, $X$ and $Y$, and all $1$-forms, $\xi$ and $\eta$,
on $A$. (Here ${\mathcal L}_X$ is the Lie derivative, {\it i.e.}, the graded
commutator $[i_X,\d]$.)
In the case of $A = TM$, the original Courant bracket \cite{Courant}
is recovered as the skew-symmetrized Dorfman bracket.

A \emph{Dirac structure} $L$ on a proto-bialgebroid is a maximally isotropic
(with respect to the canonical symmetric fiberwise bilinear form
on $A \oplus A^*$) subbundle of $A \oplus A^*$ such that $\Gamma L$ is
closed under the Dorfman or, equivalently, the Courant bracket.

\subsection{Poisson and presymplectic functions on\\ proto-bialgebroids}\label{A.2}
We shall reformulate a result of \cite{YinHe} on Dirac structures on
proto-bialgebroids, and supply an alternate, slightly more conceptual proof.
We shall utilize notations from \cite{yksPM287}.

Let $(A, \mu,\gamma, \phi, \psi)$ be a proto-bialgebroid. Let
$\pi$ (resp., $\omega$) be a bivector (resp., a
$2$-form) on $A$. Twisting the structure $\phi + \gamma +\mu +
\psi$ of $A \oplus A^*$
by $\pi$ (resp., $\omega$) yields a  proto-bialgebroid 
$\phi_\pi + \gamma_\pi +\mu_\pi + \psi_\pi$
(resp., $\phi_\omega+ \gamma_\omega +\mu_\omega +
\psi_\omega)$. In terms of the big bracket, the explicit formulas are
as follows \cite{R2002} \cite{T} \cite{yksPM287}. 

\noindent $\bullet$ For a bivector, $\pi$, of bidegree $(2,0)$,
\begin{equation}\label{phi}
\left\lbrace
\begin{array}{ll}
\vspace{.1cm}
\phi_{\pi}
= \phi - \{\gamma, \pi\} + \frac{1}{2} \{\{\mu, \pi\},\pi\} -
\frac{1}{6}  
\{\{\{\psi, \pi\},\pi\}, \pi\}  , \\
\vspace{.1cm}
\gamma_{\pi} 
= \gamma - \{\mu, \pi\} + \frac{1}{2} \{\{\psi, \pi\},\pi\}, \\
\vspace{.1cm}
\mu_{\pi} = \mu - \{\psi,\pi\} , \\
\vspace{.1cm}
\psi_{\pi}= \psi .
\end{array}
\right.
\end{equation}


\noindent $\bullet$ For a $2$-form, $\omega$, of bidegree $(0, 2)$,

\begin{equation}\label{psi}
\left\lbrace
\begin{array}{ll}
\vspace{.1cm}
\phi_{\omega}= \phi  , \\
\vspace{.1cm} 
 \gamma_{\omega} 
= \gamma - \{\phi, \omega\}  , \\
\vspace{.1cm} 
\mu_{\omega} = \mu - \{\gamma,\omega\}
+ \frac{1}{2} \{\{\phi, \omega\},\omega\}  , \\
\vspace{.1cm} 
\psi_{\omega}
= \psi - \{\mu, \omega\} + \frac{1}{2} \{\{\gamma, \omega\},\omega\} -
\frac{1}{6}  
\{\{\{\phi, \omega\},\omega\}, \omega\}  . 
 \end{array}
\right.
\end{equation}

A bivector $\pi$ (resp., a $2$-form $\omega$) is 
called a \emph{Poisson function}
(resp., a \emph{presymplectic function}) if $ \phi_\pi = 0$
(resp., $\psi_\omega =0$). A bivector is a Poisson function
if and only if the twisted proto-bialgebroid is a quasi-Lie
bialgebroid, $(A,\mu_\pi,\gamma_\pi, 0,\psi)$. A $2$-form is a 
presymplectic function if and only if the
twisted proto-bialgebroid is a 
Lie-quasi bialgebroid, $(A,\mu_\omega,\gamma_\omega, \phi,0)$.
The following proposition \cite{R2002}
\cite{T} \cite{yksPM287} extends results of \cite{LWX}.
For (i), see also Yin and He \cite{YinHe}, theorem 5.6.

\begin{proposition}\label{pfunctions}
(i) The graph of a bivector $\pi$ is a Dirac structure if and only if $\pi$ is a
Poisson function.

(ii) The graph of a $2$-form $\omega$ is a Dirac structure if and only if $\omega$ is a
presymplectic function.
\end{proposition}
 
\subsection{The characteristic pair of a Dirac structure}\label{A.3}
In \cite{YinHe}, Yin and He obtained the characterization of Dirac
structures in proto-bialgebroids in terms of characteristic pairs, thus
generalizing results of Liu \cite{Liu}.
We shall show that 
their main result (theorems 4.6 and 5.5 of \cite{YinHe}) 
and the dual statement 
can be formulated in terms of
Poisson and presymplectic functions, and we shall sketch a simple proof.

Let $D$ be a subbundle of a proto-bialgebroid,
$(A, \mu,\gamma, \phi, \psi)$, let  $D^\perp$ be its
orthogonal in $A^*$, 
and let $\pi$ be a bivector on $A$.

\begin{definition} 
We say that a bivector $\pi$ is a \emph{Poisson function
  ${\mathrm{mod}}D$} if 
$\phi_\pi \in  \Gamma(\wedge^3D)$ and $\psi_\pi \in  \Gamma(\wedge^3(D^\perp))$.
\end{definition} 

For $D= M \times \{0\} \subset A$, we recover the Poisson functions in the usual
sense \cite{T} \cite{yksPM287}, and the following theorem reduces
to (i) of Proposition~\ref{pfunctions}.

\begin{theorem}\label{proppoisson}
Let $\pi$ be a bivector on $A$. 
Let $L$ be the maximally isotropic subbundle of $A \oplus A^*$,
$$
L = \{ (X +\pi \xi, \xi)  \, | \,   X \in D, \xi \in D^\perp\}.
$$
$L$ is a Dirac structure if and only if 
$\Gamma D$ is closed under $\mu_\pi$, $\Gamma (D^\perp)$ is closed
under $\gamma_\pi$, and $\pi$ is a Poisson function ${\mathrm{mod}}D$.
\end{theorem} 

\begin{proof}
We have to find necessary and sufficient conditions for 
$$
\{\{X + \pi \xi + \xi , \phi + \gamma +\mu +\psi\}, Y + \pi \eta +
\eta \}
$$ 
for $X, Y \in \Gamma D$ and  $\xi, \eta \in \Gamma (D^\perp)$ 
to be equal to $Z + \pi \zeta + \zeta$, with  $\zeta \in \Gamma (D^\perp)$
and  $Z \in \Gamma D$.
We consider the four cases,
$\xi = \eta =0$, $\xi = Y=0$,  $X= \eta=0$,  $X = Y = 0$, i.e.,
each of the following expressions, 
\begin{equation}\label{eq1} 
\{\{X, \phi + \gamma +\mu +\psi\}, Y \},
\end{equation} 
\begin{equation}\label{othereq2} 
\{\{X, \phi + \gamma +\mu +\psi\}, \{\pi, \eta\} + \eta \}, 
\end{equation}
\begin{equation}\label{eq3}
\{\{\{\pi, \xi\} + \xi , \phi + \gamma +\mu +\psi\}, Y \},
\end{equation} 
\begin{equation}\label{eq4}
\{\{\{\pi, \xi\} + \xi , \phi + \gamma +\mu +\psi\}, \{\pi, \eta\} +
\eta \}.
\end{equation} 
We use the Jacobi identity
to write each expression as $W + \zeta$, with $W \in \Gamma A$ and
$\zeta \in \Gamma (A^*)$, and in each case we derive the condition for
$W - \pi \zeta$
to be in $\Gamma D$. For this we write that the duality bracket of $W - \pi
\zeta$ with any element $\chi \in \Gamma (D^\perp)$ vanishes.
In the case of \eqref{eq1}, we find  the condition
$\phi_\pi \in \Gamma (\wedge^3 D)$.
In the case of \eqref{othereq2}, we find the condition
$[ \Gamma (D^\perp),  \Gamma (D^\perp)]^{\gamma_\pi} \subset \Gamma (D^\perp)$.
In the case of \eqref{eq3}, we find the condition
$[ \Gamma D,  \Gamma D]^{\mu_\pi} \subset \Gamma D$.
In the case of \eqref{eq4}, we find  the condition
$\psi_\pi \in \Gamma (\wedge^3 D^\perp)$.
\end{proof}

When $L$ is a Dirac structure, $(\pi, D)$ is called the \emph{characterstic
pair} for $L$ \cite{Liu} \cite{YinHe}. (Each Dirac structure $L$ such
that $L \cap A$ has constant rank is defined by a characterisitc pair
\cite{Liu}.)
For {\it dual characterisitic pairs}, $(\omega, F)$, 
we consider a subbundle  $F$ of $A^*$ and a $2$-form $\omega$. 

\begin{definition} 
We say that a $2$-form 
$\omega$ is a \emph{presymplectic function ${\mathrm{mod}}F$} if 
$\phi_\omega \in \Gamma(\wedge^3(F^\perp))$ and $\psi_\omega \in
\Gamma(\wedge^3F)$.
\end{definition}

For $F= M \times \{0\} \subset A^*$, 
we recover the presymplectic functions in the usual
sense \cite{T} \cite{yksPM287}, and the following theorem reduces
to (ii) of Proposition \ref{pfunctions}.

\begin{theorem}\label{propsympl}
Let $\omega$ be a $2$-form on $A$.
Let $L$ be the maximally isotropic subbundle of $A \oplus A^*$,
$$
L = \{ (X, \xi + \omega X)  \, | \, X \in F^\perp, \xi \in F\}.
$$
$L$ is a Dirac structure if and only if 
$\Gamma (F)$ is closed under $\gamma_\omega$, $\Gamma (F^\perp)$ is closed
under $\mu_\omega$, and $\omega$ is a presymplectic function ${\mathrm{mod}}F$.
\end{theorem} 

\begin{proof} The proof is obtained from that of Theorem \ref{proppoisson} by
exchanging the roles of $A$ and $A^*$.
\end{proof}

\medskip

\noindent{\bf Outlook}. 
An alternate way to generalize bi-Hamiltonian structures 
is to consider pairs of Poisson functions on a proto-bialgebroid with
a suitably defined compatibility condition, such as the requirement that
their sum be a Poisson function.
Also, it would be interesting, when two Dirac structures are each 
defined in terms of a characteristic pair, to study a Dirac
pair in terms of a compatibility condition on the pair of
characteristic pairs.
As mentioned in the introduction, 
the case of Dirac structures in arbitrary Courant algebroids should
also be considered.
The compatibility of double Poisson structures in the sense of Van den
Bergh should be the subject of further work.
The search for examples and applications to problems in mechanics has
only begun.

\bigskip

\noindent{\bf Acknowledgments} I wish to thank Kirill Mackenzie, 
Volodya Rubtsov and Alan Weinstein for useful suggestions.

\end{document}